\documentclass[11pt,a4paper]{article}
\usepackage[utf8]{inputenc}
\usepackage{amsmath}
\usepackage{amsfonts}
\usepackage{amssymb}
\usepackage{authblk}
\usepackage{hyperref}
\usepackage{color}
\usepackage{mathtools}
\usepackage[ruled,vlined]{algorithm2e}
\usepackage{algorithmic}
\SetKwInput{KwInput}{Input}                
\SetKwInput{KwOutput}{Output}              

\author[1]{Rajeev Kumar\\ Society for Natural Technology Research,\\ 
Dept. of IT\&E, Governemnt of West Bengal, India.\\  
rajeev.ips@nltr }
\date{}

\providecommand{\keywords}[1]
{
  \small	
  \textbf{\textit{Keywords --- }} #1
}
\newcommand{\comment}[2]{}
\begin{document}
\title{A simple algorithm for finding square root modulo $p$}
\maketitle
\begin{abstract}
We propose a novel algorithm for finding square roots modulo $p$ in finite field $F^*_p$. Although there exists a direct formula to calculate square root of an element of field $F^*_p$, for   $p \equiv 3\pmod4$, but calculating square root for $p \equiv 1\pmod4$ is non trivial. Tonelli-Shanks algorithm remains the most widely used and probably the fastest when averaged over all primes [19]. This paper proposes a new algorithm for finding square roots in  $F^*_p$ for all odd primes, which shows improvement over existing method in practical terms although asymptotically gives the same run time as Tonelli-Shanks.\\
All algorithms, used in practice, for finding square root in finite field $F^*_p$ require knowledge of atleast one non-residue. However, there is no known  deterministic polynomial time method to find non-residue in $F^*_p$ where  $p \equiv 1\pmod4$. Although probability that an element, randomly  chosen, is a non-residue is $\frac{1}{2}$ (as we know that half the elements of $F^*_p$ are non-residues).\\
Apart from  efficient computation time, the proposed method does not necessarily require availability of non-residue and can work with \emph{`relative'} non-residue also. Such \emph{`relative'} non-residues are much easier to find ( probability $=\frac{2}{3}$) compared to finding non-residues ( probability $=\frac{1}{2}$).\\
\end{abstract}


\section{Introduction}
The first recorded reference of finding square root is found in the works of Bhaskaracharya (1150 AD) who considered the very special case $x^2 \equiv 30\pmod7$. However till present day,  the task of computing square roots in any finite field $F^*_p$ remains  a problem of considerable importance with application in well-known cryptographic mechanisms, like the quadratic sieve factorization method, point counting on elliptic curve(s)\cite{Li} etc.
The standard method of finding an element $\alpha$ in the finite field $F^*_p$  such that it is a quadratic residue,  is to check its Legendre symbol. For all quadratic residues, Legendre symbol evaluates to $1$ , i.e. $\genfrac(){}{0}{\alpha}{p} =1 $ or $\alpha^{(p-1)/2} \equiv 1 \pmod p$.
 Once we know that $\alpha$ is a quadratic residue (i.e. the square root of $\alpha$ exists), finding the square root modulo $p$ for $p \equiv 3\pmod4$ is straight forward and is given by $\sqrt {a} \equiv \pm a^{(p+1)/4}\pmod p$. 
However finding square root of a quadratic residue $\alpha$, for modulo $p \equiv 1\pmod 4$, is non trivial. There are explicit solutions for $p \equiv 5\pmod 8$  but the case for $p \equiv 1\pmod 8$ remains an open issue \cite {Uma},\cite {Sze}.\\
\\In the year 1891, Tonelli\cite {Tonelli}  proposed an algorithm to find square roots modulo $p$, followed by Cipolla's \cite {Cipolla} algorithm in 1903.  In the year 1972, Daniel shanks improved upon Tonelli's algorithm and presented an efficient Tonelli-Shanks \cite {Shanks} algorithm. This algorithm has found wide spread acceptance. It starts with approximate square root and improves it  by finding better approximation  till it finds an exact square root. Few other algorithms of Adleman-Manders-miller \cite {AMM77}, a generalisation of Tonelli-shanks for taking $r^{th}$ root  in 1977, followed by M.O.Rabin in the year 1980 (Berlekamp-Rabin \cite {Rabin}) and in the year 1986  Peralta \cite {Peralta} dealt with issue of finding square root.\\ All these algorithms are efficient but require knowledge of one non-residue for finding the square root, hence are considered probabilistic. However under the assumption that ERH is true, Ankeny \cite {Ankeny}  showed that the least quadratic non-residue over $F^*_p$ is less than $c\log^2{p}$ for some constant c. This implies that all probabilistic algorithms for finding a quadratic non-residue as mentioned  previously can be improved to a deterministic polynomial-time algorithms, if ERH is true. It is easy to see that the least quadratic non residue must be a prime \cite {Sze}.
Schoof \cite { Schoof}, in 1985, used elliptic curves to propose a deterministic algorithm to find square roots modulo  prime. This algorithm is efficient (polynomial time) for some residues but not in general. In 2011, Tsz-Wo Sze \cite {Sze}  proposed a deterministic algorithm to find square roots over finite fields without being given any quadratic non-residue.  However, run time of such an algorithm is  $\mathcal{O}(\log^3{p})$ which becomes computationally expensive and impractical for larger primes. 
\\In this paper we propose an algorithm (probabilistic - as it require a non -residue  ) to deal with square root of quadratic residue for all odd primes. In the proposed new algorithm, the traditional method of dealing with $p \equiv 3\pmod4$  has been extended for $p \equiv 1\pmod4$, with the help of use of quadratic non-residue $\beta$ of the same finite field $F^*_p$. Please see the section with 'Future Ideas' for more details. This new algorithm is much simpler, but equally efficient. Given a quadratic non-residue, proposed algorithm on average takes  $\mathcal{O}(r^2)$ time and for special primes such that $2^r>d \text{ where } p=2^rd+1$ approaches $\mathcal{O}(\log^2{p})$.
\\In fact in some cases depending upon the element whose square root is sought, algorithm works even without the availability of strict non quadratic residue.

\section{Preliminary Concepts }

The basic ideas on which this algorithm relies are - 
\subsection{ Conditions for quadratic residue in a Field}
We know that for square root of  an element $\alpha$ of   finite field $F^*_p$  to exist,  Legendre symbol  $\genfrac(){}{0}{\alpha}{p} =1 $ or $\alpha^{(p-1)/2} \equiv 1 \pmod p$. $\alpha$ is also called quadratic residue.\\
\\ 
Similarly an element $\beta$  of finite field $F^*_p$  is a non-residue if Legendre symbol  $\genfrac(){}{0}{\beta}{p} = -1 $ or $\beta^{(p-1)/2} \equiv -1 \pmod p$.\\
\subsection{ Number of roots of unity in Finite Field}
It is also a proven fact that there are only two square roots of unity in a finite field $F^*_p$  viz. $ \{1, p-1\}$.
\subsection{ Conditions for  Group $\mathbb{Z}_{p^k}^{*}$}
It can be easily seen that for multiplicative group $\mathbb{Z}_{p^k}^{*}$, roots of unity are only two viz. $\{1,p^k-1\}$. For any $ \alpha \in \mathbb{Z}_{p^k}^{*} $, $\alpha$ is a quadratic residue, i.e. $\sqrt{\alpha}$ exists if $\alpha^{\phi{(p^k)}/2} \equiv 1 \pmod p$. Similarly for any $ \beta \in \mathbb{Z}_{(p^k)}^{*} $, $\beta $ is a quadratic non-residue, i.e. $\sqrt{\beta}$ does not exist if $\beta ^{\phi{(p^k)}/2} \equiv -1 \pmod p $. Here $\phi(p^k) \text{ is euler's totient } = (p-1)p^{k-1} $.
\subsection{\emph{`relative'} non-residue}
Let $\alpha,\beta,\gamma,\gamma_1,\gamma_2 \in F^*_p $ and $ p-1= 2^rd$ where  $ r,d,j $ positive integers such that $  r,d \geq 1$ and $d$ an odd integer. $i$ is a non-negative integer. Also, let $\alpha $ be quadratic residue, $\beta $ be non-residue and $\gamma,\gamma_1,\gamma_2$ be any elements of $F^*_p $. We define another function with domain = $F_p^*$ and range $\{-1,0,1,\cdots r-1\}$\ -\\
\\
$f(\gamma) =  
\left\{
\begin{array}{ll}
-1 & \mbox{if } \gamma^d \equiv  1 \pmod p \\
i  & \mbox{if } \gamma^{2^id} \equiv -1 \pmod p\\
\end{array}
\right.$\\
\\
Let us call it f-value of element $ \gamma$. Please note this is different than order of element in finite field $F^*_p$. Also unlike order of an element, this function can be easily calculated in $\log {p} $ time. \textbf{if $f(\alpha) > f(\gamma) $, we would say $\gamma$ is \emph{`relative'} non-residue to $\alpha$ in field $F_p^*$.}

We also note the following properties of this function :- \\
\begin{itemize}
		\item $f(\beta)= r-1$.\\
		By definition $\beta $ being non-residue means $\beta^{2^{r-1}d} \equiv -1 \pmod p\text{ hence } f(\beta)=r-1 $. 
		\item $-1\leq f(\alpha) \leq r-2$\\
		By definition $\alpha $ being residue means $\alpha^{2^{r-1}d} \equiv 1 \pmod p \text{ hence } f(\alpha)<r-1 $.\\
		Hence we get following range  of values for function  $-1\leq f(\alpha) \leq r-2$.\\
		\item $f(\beta) > f(\alpha)$.\\
		As $ f(\beta)=r-1 $ and $-1\leq f(\alpha) \leq r-2$. \\
		This implies $f(\beta) > f(\alpha)$.
		\item if $f(\gamma)=0$ then $f(-\gamma)= -1$.\\
		$f(\gamma)=0$ means $\gamma^d \equiv-1 \pmod p $, this implies $(-\gamma)^d \equiv 1 \pmod p$ as $d$ is an odd integer.\\
		Hence $f(-\gamma)= -1 $  ( by function definition).
		\item if $f(\gamma)=-1$ then $f(-\gamma)= 0$.\\
		$f(\gamma)=-1$ means $\gamma^d \equiv 1 \pmod p$, this implies $(-\gamma)^d \equiv -1 \pmod p$ as $d$ is an odd integer.\\
		Hence $f(-\gamma)= $0  ( by function definition).
		\item $f(\gamma) = f(-\gamma) $ for $f(\gamma)>0$.\\
		Let $f(\gamma)= i$ and $i>0$, that would mean $ \gamma^{2^id}\equiv -1 \pmod p$.\\
		This implies $f(-\gamma)= i$ as $\gamma^{2^id}=(-\gamma)^{2^id} $ as $i>0$.\\
		Hence $f(\gamma) = f(-\gamma)$ for $f(\gamma)>0$.
		\item $f(\gamma) = f(\gamma^{-1}) $.\\
		For $f(\gamma) = -1$ implies $\gamma^d \equiv 1 $\\
		Now say ${(\gamma^{-1})}^d \equiv \phi $ implies $\phi.\gamma^d \equiv 1$ implies $\phi=1$.\\
		For non trivial cases, $f(\gamma) = \delta$ where $\delta \geq 0 $ implies $\gamma^{2^\delta.d }\equiv -1 $ \\
		Now say ${(\gamma^{-1})}^{2^\delta.d }\equiv \phi $ implies $\phi.\gamma^{2^\delta.d } \equiv 1$ implies $\phi=-1$.\\
		Hence $f(\gamma) = f(\gamma^{-1}) $
		\item if $f(\alpha)=-1$ then $\sqrt\alpha \equiv {\alpha}^{(d+1)/2} \pmod p $.\\
		$f(\alpha)=-1$ means $\alpha^d \equiv 1 \pmod p$ implies $\alpha^{d+1} \equiv \alpha \pmod p.\\ $
		Hence $\sqrt\alpha \equiv {\alpha}^{(d+1)/2} \pmod p $.
		\item $f(\sqrt{\alpha}) =f(\alpha)+1 $ for $ f(\alpha) \geq 0$.\\
		Let us suppose  $f(\sqrt{\alpha})=j$ that would mean ${(\sqrt{\alpha})}^{2^jd} \equiv -1 \pmod p.$\\
		For  $j\geq 1 $, this would mean  ${(\alpha)}^{2^{j-1}d} \equiv -1$. \\
		Now if  $f(\alpha)=i$ and  $i\geq 0$ means $\alpha^{2^id} \equiv -1 \pmod p.$\\
		Comparing two expressions implies  $j-1=i $, i.e. $ j=i+1$ for $i \geq 0 $.\\
		Hence $f(\sqrt{\alpha}) =f(\alpha)+1 $ for $ f(\alpha) \geq 0$. 
		
		\item $f(\gamma_1\gamma_2) < f(\gamma_1) $   if $f(\gamma_1)=f(\gamma_2) $ for $f(\gamma_1),f(\gamma_2) \geq 0 $.\\
		Let us take $f(\gamma_1)=f(\gamma_2)=i$ and $ i\geq 0$. This implies $\gamma_1^{2^id} \equiv \gamma_2^{2^id} \equiv -1 \pmod p$.\\
		This means $(\gamma_1\gamma_2)^{2^id} \equiv 1 \pmod p$. Hence $f(\gamma_1\gamma_2) <i$.\\
		Hence $f(\gamma_1\gamma_2) < f(\gamma_1) $ if $f(\gamma_1)=f(\gamma_2)\geq0$.\\
		If $f(\gamma_1)=-1$, nothing can be said about value of $f(\gamma_1\gamma_2)$ in relation to value of $f(\gamma_1),f(\gamma_2)$.
		\item $f(\gamma_1\gamma_2) = \text{max. of }(f(\gamma_1),f(\gamma_2)) $ if $f(\gamma_1) \neq f(\gamma_2) $.\\
		Let $f(\gamma_1) =i \text{ and } f(\gamma_2)= j \text{ where } ,i,j\geq 0 \text{ and } i\neq j$.\\
		Without loss of generality we take $ i >j$.
		Now $\gamma_1^{2^id} \equiv -1 \pmod p$ and $\gamma_2^{2^jd} \equiv -1 \pmod p$.\\
		Please note as $i>j$ this implies  $\gamma_2^{2^id} \equiv 1$ \\
		This implies ${(\gamma_1\gamma_2)}^{2^{i}d} \equiv \gamma_1^{2^id}\gamma_2^{2^id} \equiv -1.1  \equiv -1$.\\
		Hence $f(\gamma_1\gamma_2) =i=\text{ maximum of }(f(\gamma_1),f(\gamma_2)). $\\	
\end{itemize}
\comment{
Let  $\alpha_1,\alpha_2,p_i^{e_i} \in F^*_p$. Now, it is easy to see, from above mentioned assertions that for a given $ \alpha_1$, the only way we can arrive at  $ \alpha_2$ such that $f(\alpha_2) > f(\alpha_1)$, is by  -\\
\begin{itemize}
\item[a] Taking square root as $f(\sqrt{\alpha_1}) =f(\alpha_1)+1 $ for $ f(\alpha_1) \geq 0$.
\item[b] Finding prime factor as $f(p_i^{e_i}) < f(p_i)$.
\item[c] Finding factors with equal function value as $f(\alpha_1\alpha_2) < f(\alpha_1) $   if $f(\alpha_1)=f(\alpha_2) $
\end{itemize}
}

The possible values of $ f(\gamma) \in S_v=\{-1,0,1,2,3,\cdots,r-2,r-1 \}$(from definition of function).
Let $S_{f_i} $ denote the set all elements in $F_p^*$ having f-value equal to $i$, i.e.\\$S_{f_i}= \{ x: x \in Z_p^* \text{ such that } f(x)=i \} $. It is easy to see that if $ x \in S_{f_i} $ then $ \pm \sqrt{x} \in S_{f_{i+1}}$ for $0 \leq i \leq r-2 $. And if $ x \in S_{f_{-1}} $ then $\sqrt{x} \in S_{f_{-1}} $ and  $-\sqrt{x} \in S_{f_{0}} $. Using these facts it is easy to calculate the cardinality of set $S_{f_i} $ as\\
\\
$|S_{f_i}| =  
\left\{
\begin{array}{ll}
2^{i}d & \mbox{if } i\geq 0 \\
i  & \mbox{if } i=-1\\
\end{array}
\right.$\\
\\
Please note that $\sum_{i=-1}^{r-1} |S_{f_i}| = \sum \left( d+d+2d+2^2d+2^3d+\cdots 2^{r-1}d \right) = d+ d(2^0+2^1+2^2+\cdots 2^{r-1})= d+d(2^r-1) = 2^rd = p-1$, equal to total elements in field $F_p^*$.\\

\subsection{Theorem 1: for quadratic residue modulo $p$}
Let  $p $ be an odd prime. Let $p-1$ be expressed as $p-1=2^rd$ where d  odd integer ($d\geq 1 $)  and $r$  be any integer ($r >0) $. Let $\alpha$  be a quadratic residue and let $\beta$  be quadratic non residue in the finite field  $\mathbb{F}_p^*$. \textbf{Then a non-negative integer $m$ can always be found such that $\alpha^d\beta^{2md} \equiv 1 \pmod p $}. This will imply $\alpha^{(d+1)}\beta^{2md} \equiv \alpha \pmod p $ and hence $  \sqrt{\alpha} \equiv \pm (\alpha)^{(d+1)/2}\beta^{md} \pmod p $.\\
Let  $p \in \mathbb{P}$ be an odd prime. Let $p-1$ be expressed as $p-1=2^rd$ where $d$  odd integer ($d\geq 1 $)  and $r$  be any integer $(r >0) $. Let $\alpha$  be a quadratic residue and let $\beta$  be quadratic non residue in the finite field  $\mathbb{F}_p^*$. \textbf{Then a non-negative integer $m$ can always be found such that $\alpha^d\beta^{2md} \equiv 1 \pmod p $}. This will imply $\alpha^{(d+1)}\beta^{2md} \equiv \alpha \pmod p $ and hence $  \sqrt{\alpha} \equiv \pm (\alpha)^{(d+1)/2}\beta^{md} \pmod p $.\\
\\
\textbf{Proof :} \\
Given a non-quadratic residue $\beta $ and a quadratic residue $\alpha $ such that  $\alpha^{2^{r-1}d} \equiv 1 $. 
We  will reduce  $\alpha ^{2^{r-1}d}$  to $\alpha ^{d}$ by successive square rooting operation in $(r-1)$ steps. 
Starting from $\alpha^{2^{r-1}d} \equiv 1 $, taking square root gives $\alpha^{2^{r-2}d} \equiv \pm 1 $. At each step we multiply expression with $\beta^{\lambda_i2^{r-1}d} $, where $\lambda_i \in \{0,1\} $. we choose value of $\lambda_i=0$ if square root was 1 and we choose $\lambda_i=1$ if the square root was $-1$.\\
We note that $\beta^{2^{r-1}d} \equiv -1$ as $\beta$ is non residue. 
Hence after $k$ steps we get the following expression where $\lambda_i \in \{0,1\} $ :- 

\begin{align*}
&\alpha^d*\beta^{\lambda_12^{1}d}*\beta^{\lambda_22^{2}d} \cdots \beta^{\lambda_{r-1}2^{r-1}d} &\equiv 1 \pmod p\\ 
&\alpha^d*\beta^{(\lambda_12^{1}d +\lambda_22^{2}d+ \cdots \lambda_{r-1}2^{r-1}d)} &\equiv 1 \pmod p\\  
&\alpha^d*\beta^{2d(\lambda_1 +\lambda_22^1+\cdots \lambda_{r-1}2^{r-2})} &\equiv 1 \pmod p\\ 
&\alpha^d*\beta^{2md} &\equiv 1 \pmod p
\end{align*}
where  $m = (\lambda_1 +\lambda_22^1+ \lambda_{r-1}2^{r-2})$  and $\lambda_i \in \{0,1\} $. Also note that  for $r=1$ theorem is true as $ m=0$. This proves the theorem.
\subsection{Theorem 2: Set of possible solutions mod $p$}
Let  $p \in \mathbb{P}$ be an odd prime, such that  $p-1=2^rd$ where $d$ be an odd integer ($d\geq 1 $)  and $r$  be any integer ($r >0) $. Let $\alpha$  be a quadratic residue and let $\beta$  be quadratic non-residue in the finite field $F^*_p$. Then  $  \sqrt{\alpha} \pmod p $ can always be found in the set  $\{ \alpha^{(d+1)/2}\beta^{2^{(r-1-k)}id}  \text{ for } 0 \leq i \leq 2^k-1  \}$  where  $k$, for which $\alpha^{2^{k}d} \equiv 1 $ and $0\leq k \leq (r-1)$. \\ 
\\\textbf{Proof :} \\
Given a non-quadratic residue $\beta $ and a quadratic residue $\alpha $ such that  $\alpha^{2^{k}d} \equiv 1 $ for  $0\leq k \leq (r-1)$. 
We  will reduce  $\alpha ^{2^{k}d}$  to $\alpha ^{d}$ by successive square rooting operation in $k$ steps. 
Starting from $\alpha^{2^{k}d} \equiv 1\pmod p $, taking square root gives $\alpha^{2^{k-1}d} \equiv \pm 1 \pmod p $. At each step we multiply expression with $\beta^{\lambda_i2^{r-1}d} $, where $\lambda_i \in \{0,1\} $. we choose value of $\lambda_i=0$ if square root was 1 and we choose $\lambda_i=1$ if the square root was $-1$.\\
We note that $\beta^{2^{r-1}d} \equiv -1$ as $\beta$ is non residue. 
Hence after $k$ steps we get the following expression where $\lambda_i \in \{0,1\}:-$
\begin{align*}
\alpha^d*\beta^{\lambda_12^{r-k}d}*\beta^{\lambda_22^{r-k+1}d}* \cdots * \beta^{\lambda_k2^{r-1}d} \equiv 1 \pmod p\\
\alpha^{d+1}*\beta^{\lambda_12^{r-k}d}*\beta^{\lambda_22^{r-k+1}d} * \cdots * \beta^{\lambda_k2^{r-1}d} \equiv \alpha \pmod p
\end{align*}
Rearranging terms
\begin{align*}
\sqrt{\alpha } \pmod p&\equiv \alpha^{(d+1)/2}*\beta^{\lambda_12^{r-k-1}d}*\beta^{\lambda_22^{r-k}d}* \cdots *\beta^{\lambda_k2^{r-2}d} \\
\sqrt{\alpha } \pmod p &\equiv \alpha^{(d+1)/2}*\beta^{\lambda_k2^{r-2}d}*\beta^{\lambda_{k-1}2^{r-3}d} *\cdots * \beta^{\lambda_12^{r-k-1}d} 
\end{align*}
It is easy to see that possible values of $\sqrt{\alpha } $ can be all possible permutations of $\lambda_i $. As there are only two possible values of $\lambda_i \in \{0,1\} $, hence value of $\sqrt{\alpha }$ has to be equal to  one of these possible $2^k $ values. We also note that all other terms of $\beta $ evaluate to $1$(i.e. exponent $=0$, when $ \lambda_i=0$) or can be obtained by successive squaring of $\beta^{2^{r-k-1}d}$. Basically it is $k$-bit number, where value of bit is either 1 (exponent of $ \beta =0$) or $(\beta^{2^{r-k-1}d})^{2^j}$ ($j^{th}$ place value ) for the $j^{th}$ placed bit, where $0\leq j\leq k $. Hence all possible solutions are   $\beta^{2^{(r-1-k)}id} $ for $0\leq i \leq 2^k-1$.
Hence set of all possible values of $\sqrt{\alpha } $  is equal to $\{ \alpha^{(d+1)/2}\beta^{2^{(r-1-k)}id} \text{ for } 0 \leq i \leq 2^k-1 \} $.  We also note that as $\alpha  $ is a quadratic residue, hence there exists a  $0 \leq k \leq (r-1) $ such that  $\alpha^{2^kd \equiv 1 } $. This proves the theorem.\\ 
\subsection{Theorem 3: for quadratic residue modulo $p^k$}
Let  $p \in \mathbb{P}$ be an odd prime and let $k$ be an integer $k>0$. Let $\phi(p^k)=2^RD$ where $D$ be an  odd integer ($D\geq 1 $)  and $R$ be any integer ($R >0) $. Let $\alpha$  be a quadratic residue and let $\beta$  be quadratic non-residue in the multiplicative group $\mathbb{Z}_{p^k}^{*}$. \textbf{Then a non-negative integer $m\geq0$ can always be found such that $\alpha^D\beta^{2mD} \equiv 1  \pmod {p^k} $}. This will imply $\alpha^{(D+1)}\beta^{2mD} \equiv \alpha \pmod {p^k} $ and hence $  \sqrt{\alpha} \equiv \pm (\alpha)^{(D+1)/2}\beta^{mD} \pmod {p^k} $.\\
\\
\textbf{Proof :} \\
Proof is similar to Theorem 1 by noting that roots of unity in multiplicative group $\mathbb{Z}_{p^k}^{*}$ are $ \{1,-1\}$. Also $\alpha$ is a quadratic residue hence  $\alpha^{\phi{(p^k)}/2} \equiv 1\pmod {p^k} $. Similarly for any non residue  $\beta $  if $\beta ^{\phi{(p^k)}/2} \equiv -1 \pmod {p^k} $. Here $\phi(p^k) \text{ is euler's totient } = p^{(k-1)}(p-1) $. Arguing on similar lines as Theorem 1, it eventually leads us to $\alpha^D\beta^{2mD} \equiv 1  \pmod {p^k} $ where $ m\geq 0$.

\subsection{Theorem 4: Set of possible solutions mod $p^k$}
Let  $p \in \mathbb{P}$ be an odd prime and let $k$ be an integer $k>0$. Let $\phi(p^k)=2^RD$ where $D$ be an  odd integer ($D\geq 1 $)  and $R$ be any integer ($R >0) $. Let $\alpha$  be a quadratic residue and let $\beta$  be quadratic non residue in the multiplicative group $\mathbb{Z}_{p^k}^{*}$. Then  $\sqrt{\alpha} \pmod {p^k} $ can always be found in the set  $\{ (\alpha)^{(D+1)/2}\beta^{2^{(R-1-j)}iD} \text{ for } 0 \leq i \leq 2^j-1 \}$ where  $j$ is such that  for which $\alpha^{2^{j}D} \equiv 1$ and $0\leq j \leq (R-1)$  .  \\ 
\\
\textbf{Proof :} \\
This can be proven on the same lines as  Theorem 2. Here we  try to prove it by noting that here $\phi(p^k)=2^RD$ and also as given value of $j$ for which $\alpha^{2^{R-1-j}D} \equiv 1$.  
Arguing on similar lines as Theorem 2, we get  $  \sqrt{\alpha} \pmod{p^k} $ can always be found in the set  $\{ \alpha^{(D+1)/2}\beta^{2^{(R-1-j)}iD} \text{ for } 0 \leq i \leq 2^j-1 \}$ where  $j$ is such that  for which $\alpha^{2^{j}D} \equiv 1$  and $0\leq j \leq (R-1)$.
\subsection{Theorem 5: \emph{`relative'} non-residue is sufficient to calculate square root}
Let $\alpha,\gamma \in F^*_p $ and $ p-1= 2^rd$ where  $  r,d \geq 1$ and $d$ an odd integer. Also, let $\alpha $ be quadratic residue  and $\gamma$ be any element of $F^*_p $. If  $f(\gamma) > f(\alpha) $ (function f as defined under \emph{`relative'} non-residue section), we can always find $m \geq 0$ such that $\alpha^d\gamma^{2md} \equiv 1 $ hence $\sqrt{\alpha}\equiv \alpha^{d+1/2}\gamma^{md}$. \\ 
\\
\textbf{Proof :} \\
Let $f(\alpha) = s $ and $f(\gamma) = t $ . Given that $ t>s $. Please refer to definition of function f, given in section 2.4 \emph{`relative'} non-residue above, we note that $f(\alpha) = s $ implies  $\alpha^{2^sd} \equiv -1$. Also note $ \gamma^{2^td} \equiv -1$ and $ t > s$. We reduce  $\alpha ^{2^{r-1}d} \equiv 1 $ (as $\alpha$ is quadratic residue)  to $\alpha ^{d}$ by successive square rooting operation in $(r-1)$ steps. At each step we multiply expression with $\gamma^{\lambda_i2^{j-1}d} $, where $\lambda_i \in \{0,1\} $. We choose value of $\lambda_i=0$ if square root is 1 and we choose $\lambda_i=1$ if square root is $-1$. However, $ f(\alpha)=s$ implies $\alpha^{2^{s+1}d} \equiv 1 $, which implies $\lambda_i=0 $ for $s <i  $ \\
Hence after $(r-1)$ steps we get \\
$\alpha^d*\beta^{\lambda_12^{1}d}*\beta^{\lambda_22^{2}d} \cdots \beta^{\lambda_{r-1}2^{r-1}d} \equiv 1 $.

 Hence we get $ \alpha^{2^id}\gamma^{2^jd} \equiv 1$. As $j>i$ , we can argue exactly on same lines as done in proof of Theorem 1 (section 2.5) to show that we can always find  $m>0$ such that $\alpha^d\gamma^{2md} \equiv 1 $ hence $\sqrt{\alpha}\equiv \alpha^{d+1/2}\gamma^{md}$. 
Given a non quadratic residue $\beta $ and a quadratic residue $\alpha $ such that  $\alpha^{2^{r-1}d} \equiv 1 $. 

Starting from $\alpha^{2^{r-1}d} \equiv 1 $, square root will give $\alpha^{2^{r-2}d} \equiv \pm 1 $. 
We note that $\beta^{2^{r-1}d} \equiv -1$ as $\beta$ is non residue. 
hence after k steps we get the following expression where $\lambda_i \in \{0,1\} $ :- \\
$\alpha^d*\beta^{\lambda_12^{1}d}*\beta^{\lambda_22^{2}d} \cdots \beta^{\lambda_{r-1}2^{r-1}d} \equiv 1 $.
\\ Implying, 
$\alpha^d*\beta^{(\lambda_12^{1}d +\lambda_22^{2}d+ \cdots \lambda_{r-1}2^{r-1}d)} \equiv 1 $.
 \\ Implying, 
$\alpha^d*\beta^{2d(\lambda_1 +\lambda_22^1+ \cdots \lambda_{r-1}2^{r-2})} \equiv 1 $.
  \\ Implying, 
$\alpha^d*\beta^{2md} \equiv 1 $ where  $m = (\lambda_1 +\lambda_22^1+ \cdots \lambda_{r-1}2^{r-2})$  and $\lambda_i \in \{0,1\} $. Also note that  for $r=1$ theorem is true as $ m=0$. This proves the theorem.
\comment{
We also note the following properties of this function- \\
\begin{itemize}
		\item $f(\beta)= r-1$.\\
		By definition $\beta $ being non-residue means $\beta^{2^{r-1}d} \equiv -1 \text{ hence } f(\beta)=r-1 $. 
		\item $-1\leq f(\alpha) \leq r-2$\\
		By definition $\alpha $ being residue means $\alpha^{2^{r-1}d} \equiv 1 \text{ hence } f(\alpha)<r-1 $.\\
		Hence we get following range  of values for function  $-1\leq f(\alpha) \leq r-2$.\\
		\item $f(\beta) > f(\alpha)$.\\
		As $ f(\beta)=r-1 $ and $-1\leq f(\alpha) \leq r-2$. \\
		This implies $f(\beta) > f(\alpha)$.
		\item if $f(\gamma)=0$ then $f(-\gamma)= -1$.\\
		$f(\gamma)=0$ means $\gamma^d \equiv=-1 $, this implies $(-\gamma)^d \equiv 1 $ as $d$ is an odd integer.\\
		Hence $f(-\gamma)= -1 $  ( by function definition).
		\item if $f(\gamma)=-1$ then $f(-\gamma)= 0$.\\
		$f(\gamma)=-1$ means $\gamma^d \equiv= 1 $, this implies $(-\gamma)^d \equiv -1 $ as $d$ is an odd integer.\\
		Hence $f(-\gamma)= $0  ( by function definition).
		\item $f(\gamma) = f(-\gamma) $ for $f(\gamma)>0$.\\
		Let $f(\gamma)= i$ and $i>0$, that would mean $ \gamma^{2^id}\equiv -1$.\\
		This implies $f(-\gamma)= i$ as $\gamma^{2^id}=(-\gamma)^{2^id} $ as $i>0$.\\
		Hence $f(\gamma) = f(-\gamma)$ for $f(\gamma)>0$.
		\item if $f(\alpha)=-1$ then $\sqrt\alpha \equiv {\alpha}^{(d+1)/2} $.\\
		$f(\alpha)=-1$ means $\alpha^d \equiv 1$ implies $\alpha^{d+1} \equiv \alpha.\\ $
		Hence $\sqrt\alpha \equiv {\alpha}^{(d+1)/2} $.
		\item $f(\sqrt{\alpha}) =f(\alpha)+1 $ for $ f(\alpha) \geq 0$.\\
		Let us suppose  $f(\sqrt{\alpha})=j$ that would mean ${(\sqrt{\alpha})}^{2^jd} \equiv -1.$\\
		For  $j\geq 1 $, this would mean  ${(\alpha)}^{2^{j-1}d} \equiv -1$. \\
		Now if  $f(\alpha)=i$ and  $i\geq 0$ means $\alpha^{2^id} \equiv -1.$\\
		Comparing two expressions implies  $j-1=i $, i.e. $ j=i+1$ for $i \geq 0 $.\\
		Hence $f(\sqrt{\alpha}) =f(\alpha)+1 $ for $ f(\alpha) \geq 0$. 
		
		\item $f(\gamma_1\gamma_2) < f(\gamma_1) $   if $f(\gamma_1)=f(\gamma_2) $ for $f(\gamma_1),f(\gamma_2) \geq 0 $.\\
		Let us take $f(\gamma_1)=f(\gamma_2)=i$ and $ i\geq 0$. This implies $\gamma_1^{2^id} \equiv \gamma_2^{2^id} \equiv -1 $.\\
		This means $(\gamma_1\gamma_2)^{2^id} \equiv 1$. Hence $f(\gamma_1\gamma_2) <i$.\\
		Hence $f(\gamma_1\gamma_2) < f(\gamma_1) $ if $f(\gamma_1)=f(\gamma_2)\geq0$.\\
		If $f(\gamma_1)=-1$, nothing can be said about value of $f(\gamma_1\gamma_2)$ in relation to value of $f(\gamma_1),f(\gamma_2)$.
		\item $f(\gamma_1\gamma_2) = \text{max. of }(f(\gamma_1),f(\gamma_2)) $ if $f(\gamma_1) \neq f(\gamma_2) $.\\
		Let $f(\gamma_1) =i \text{ and } f(\gamma_2)= j \text{ where } ,i,j\geq 0 \text{ and } i\neq j$.\\
		Without loss of generality we take $ i >j$.
		Now $\gamma_1^{2^id} \equiv -1$ and $\gamma_2^{2^jd} \equiv -1$.\\
		Please note as $i>j$ this implies  $\gamma_2^{2^id} \equiv 1$ \\
		This implies ${(\gamma_1\gamma_2)}^{2^{i}d} \equiv \gamma_1^{2^id}\gamma_2^{2^id} \equiv -1.1  \equiv -1$.\\
		Hence $f(\gamma_1\gamma_2) =i=\text{ maximum of }(f(\gamma_1),f(\gamma_2)). $\\	
\end{itemize}

Let  $\alpha_1,\alpha_2,p_i^{e_i} \in F^*_p$. Now, it is easy to see, from above mentioned assertions that for a given $ \alpha_1$, the only way we can arrive at  $ \alpha_2$ such that $f(\alpha_2) > f(\alpha_1)$, is by  -\\
\begin{itemize}
\item[a] Taking square root as $f(\sqrt{\alpha_1}) =f(\alpha_1)+1 $ for $ f(\alpha_1) \geq 0$.
\item[b] Finding prime factor as $f(p_i^{e_i}) < f(p_i)$.
\item[c] Finding factors with equal function value as $f(\alpha_1\alpha_2) < f(\alpha_1) $   if $f(\alpha_1)=f(\alpha_2) $
\end{itemize}
}

\comment{
Let us define set $S_c = \{\text{ count of all elements for which }f(x)=i : -1 \leq i \leq (r-1)  \text{ and } x \in F_p^*  \} $.
So $S_c=\{d,d,2d,2^{2}d,\cdots 2^{r-2}d,2^{r-1}d \} $.There is one to one and onto relationship between two sets $S_c,S_v$. Combining both sets we get the following set $ S=$ \{($x \in S_v$, $y \in S_c$)\} -\\
$S=\{(-1,d),(0,d),(1,2d),(2,2^2d),(3,2^3d),\cdots (r-2,2^{r-2}d),(r-1,2^{r-1}d)\} $.\\
This is just to re-emphasise the relative correspondence between elements of sets $S_c$ and $S_v $.\\} 
\\

\comment{
\\Also please note that set $S_{f_k}$ contains similar elements i.e. let $ Y= \{x^{2n+1} : x \in S_{f_k} , n\in \mathbb{I}\}$, then $ Y \subset S_{f_k}$. \comment{ In most of the cases it turns out that $ Y=S_{f_k}$ ( a proper analysis might reveal some hidden properties ).} Proof for $ Y \subset S_{f_k} $, is easy to see (if $ x^{2^kd} \equiv -1 $ this implies $ (x^{2n+1})^{2^kd} \equiv -1 $ hence $ x^{2n+1} \in S_{f_k}$. ), however proving that when $ Y = S_{f_k} $ is easy to see from group theory - which indicates $ Y = S_{f_k} $ whenever $gcd(2n+1,p-1)=gcd(2n+1,2^rd) \equiv 1$ . Also it is easy to note that if $x \in S_{f_k} $ then $x^{2n} \in S_{f_{k-1}}$ for $0 \leq k \leq r-1 $. 
} 

\subsection{Theorem 6: \emph{`relative'} non-residue is `easier' to find than non-residue}
Let $\alpha,\gamma \in F^*_p $ and $ p-1= 2^rd$ where  $  r,d \geq 1$ and $d$ an odd integer. Also, let $\alpha $ be quadratic residue  and $\gamma$ be any element of $F^*_p $. Given a  quadratic residue $\alpha  $  probability that a randomly chosen element $\gamma \in F^*_p $  is a \emph{`relative'} non-reside to $\alpha$ is approx. averaged over all quadratic residues.  \\ 
\\
\textbf{Proof :} \\
From the discussion under section \emph{`relative'} non-residue and using the same symbols, we note the possible values of $ f(\gamma) \in S_v=\{-1,0,1,2,3,\cdots ,r-2,r-1 \}$ (from definition of the function).\\
Let $S_{f_i} $ denote the set all elements in $F_p^*$ having f-value $=i$, i.e.\\$S_{f_i}= \{ x: x \in Z_p^* \text{ such that } f(x)=i \} $. Also  it is easy to see that if $ x \in S_{f_i} $ then $ \pm \sqrt{x} \in S_{f_{i+1}}$ for $0 \leq i \leq r-2 $. And if $ x \in S_{f_{-1}} $ then $\sqrt{x} \in S_{f_{-1}} $ and  $-\sqrt{x} \in S_{f_{0}} $. Using these facts it is easy to calculate the cardinality of set $S_{f_i} $.\\
\\
$|S_{f_i}| =  
\left\{
\begin{array}{ll}
2^{i}d & \mbox{if } i\geq 0 \\
i  & \mbox{if } i=-1\\
\end{array}
\right.$\\
\\
Please note that $\sum_{i=-1}^{r-1} |S_{f_i}| = \sum \left( d+d+2d+2^2d+2^3d+\cdots 2^{r-1}d \right) = d+ d(2^0+2^1+2^2+\cdots 2^{r-1})= d+d(2^r-1) = 2^rd = p-1$, equal to the total elements in field $F_p^*$.\\
it is clear that probability of a random element having f-value equal to $k $ is  $|S_{f_k}|/2^rd$ (=$1/2^{r-k} $ for $k \geq 0 $ and =$1/2^r $ for $k=-1$).\\ 
Also probability of getting a \emph{`relative'} non-residue, given an element $\alpha $ (with f-value$=k$) = (number of  elements having f-value $> k$) / (all elements of $F_p^* $) = $\sum_{i=k+1}^{r-1} |S_{f_i}|/2^rd $, where  $-1 \leq f(\alpha)=k \leq r-2 $. \\
Now given a random $\alpha$ (a quadratic residue , whose square root is sought, picking up a random $\gamma \in F_p^*$, such that $f(\gamma) >f(\alpha) = k$(say ), averaging it over all the possible values of k$(-1 \leq k \leq r-2)$, the required probability that $\gamma $ is suitable  is \\ 
\begin{align*}
&=1/(2^{r-1}d)(\sum_{k=-1}^{r-2}\sum_{i=k+1}^{r-1} |S_{f_k}||S_{f_i}|/2^rd )\\ 
&= 1/(2^{r-1}d)((2^rd-d)d/2^rd+\sum_{k=0}^{r-2}\sum_{i=k+1}^{r-1} 2^kd|S_{f_i}|/2^rd )\\
&= 1/(2^{r-1}d)((2^rd-d)d/2^rd+\sum_{k=0}^{r-2} (2^rd-2^{k+1}d)2^kd/2^rd)\\
&= 1/(2^{r-1})((2^r-1)/2^r+\sum_{k=0}^{r-2}(2^{r}-2^{k+1})/2^r) \\
&=1/(2^{r-1})((2^r-1)/2^r +\sum_{k=0}^{r-2}(2^r-2^{k+1})2^k/2^r)\\
&=\frac{2}{3}(1-1/2^{2r})
\end{align*}
 where minimum possible value of $r=2$ for primes $1\pmod 4$. Hence probability of finding suitable $\gamma$, is $0.625$ ( for case $r=2$) and approaches $\frac{2}{3} $ for $r>2$ (please note term $1/2^{2r}$ which approaches zero for $r>2$).
 
It is clear that probability of finding a suitable $\gamma$, given an $\alpha$, such that $f(\gamma) > f(\alpha) $ is much higher than $\frac{1}{2}$ (which is the case when we are searching for strict non residue).\\

\section{Proposed Algorithm }

Once we have determined that $\alpha$ is indeed a quadratic residue, 
the core idea is to keep taking square root starting from $\alpha^{(p-1)/2} \equiv 1 $. However when ever we get $-1$ as square root we multiply both sides by $\beta^{(p-1)/2} \equiv -1$, where $\beta$ is non residue. and continue the process again till we reach $\alpha^d*\beta^{2d}\cdots  \equiv 1 $  or $\alpha^{(2k+2)}*\beta^{2d}\cdots  \equiv \alpha $ (as $d=2k+1 $), giving us $ \alpha^{(k+1)}*\beta^{d}*\beta^{2d}\cdots \equiv \sqrt{\alpha}  $. It can be easily seen we would be able to reduce the exponent of $\alpha$ to an odd number with exponent of $\beta$ still being even. Once we achieve that, we multiply both sides by $\alpha$ and obtain an expression which has even exponent of both $\alpha$ and $\beta$ equal to $\alpha \pmod p $. This equation directly provides $\sqrt{\alpha}$. \\
\\Please note that unlike previous approaches \cite{Shanks}, \cite{AMM77}, new proposed method is basically top down traversal of binary tree. It starts from $\alpha^{(p-1)/2} \equiv 1$ and works its way down halving the exponent of $\alpha $ at each step till it reaches $\alpha^d\beta^{2md} \equiv 1 $. This is achieved in just $r-1$ steps where $p-1=2^rd$. \\
\\Same methodology can be applied, so that it works for taking square root of quadratic residues modulo $p^k$. Please see example -'E' later under the examples.\\
\subsection{Insight into the algorithm }
The core of the algorithm is based on the fact that there are only two possible square roots of any element $\alpha \in Z_p^*$ and square root of $ 1 \in \{1,p-1\}$. This allows us to construct a set with at most $2^{(r-1)}$  elements, where one of the elements must be the desired square root. Knowledge of the possible solution set and its elements available at the leaf nodes of binary tree of height at most $(r-1)$ makes top down traversal of such tree natural and fast.  \\

Runtime for this algorithm is $\mathcal{O}(r^2)$ which is asymptotically as good as Tonelli-Shanks, but ability to move up and down the search tree makes it more efficient to implement than Tonelli-Shanks method. 
\comment{
Some more saving in run time can be achieved at the expense of memory requirement , by  storing  values of $\{ \alpha^d, \alpha^{2d}, \cdots  \alpha^{2^{r-2}d} \} $  and $\{ \beta^d, \beta^{2d}, \cdots  \beta^{2^{r-2}d} \} $ .i.e memory storage requirement of $ \mathcal{O}(r)$ or in worst case, when $r$ approaches $\log{p}$, $ \mathcal{O}(\log{p})$. Run time for these special primes $p=2^r+1 $ approaches $\mathcal{O}(\log^2{p})$ in worst case. }

\subsection{Pseudo Code for square root of a quadratic residue modulo $p$}
Please note pseudo code is illustrative and no extra optimisations have been used.
\begin{algorithm}[H]
\caption{Pseudo Code  to calculate $\sqrt{\alpha} \pmod p$. Returns 0 if none exits.}
\KwInput{$\alpha,\beta \in F_p^* , p$; $\beta $ is non-residue and $ p $ a prime.}
\KwOutput{$\sqrt{\alpha} $, if it exists else  0. }

\begin{algorithmic}
\IF {$(p=2)$}
	\RETURN {$\alpha$}
\ENDIF
\IF {$(\frac{\alpha}{p})= -1$} 
	\RETURN {$0$}
\ENDIF
\STATE $\alpha pow \gets (p-1)$ and $\beta pow \gets 0$ 
\WHILE{$(\alpha pow \text{ is even })$} 
	\STATE $\alpha pow \gets \alpha pow/2 $ and $\beta pow \gets \beta pow/2$ 
	\IF {$\alpha^{\alpha pow }\beta^{\beta pow}= -1$}
		\STATE $\beta pow \gets \beta pow +(p-1)/2 $
	\ENDIF
\ENDWHILE
\STATE $\alpha pow \gets (\alpha pow +1)/2 $ and $\beta pow \gets \beta pow/2$ 
\RETURN{$\alpha^{\alpha pow }\beta^{\beta pow}$}
\end{algorithmic}
\end{algorithm}
Please note that irrespective of value of $\alpha ,\beta$ `While Loop' is executed exactly $r-1$ times. Essentially, total number of calculations inside the loop are same (if $\alpha^{\alpha pow }\beta^{\beta pow}= -1$,  one extra addition per iteration takes place.\\Inside `While loop' only basic operations (multiplications, division, addition) etc are performed and there are  no nested loops, hence run time for all ($\alpha ,\beta$) pairs remains essentially the same.\\
\comment{
Problem of finding square root of $ \alpha$ is basically a search of  in the solution set  $\{ \alpha^{(d+1)/2}\beta^{2^{(r-1)}id}  \text{ for } 0 \leq i \leq 2^{r-1}-1  \}$ whose cardinality is $ 2^{r-1}$. \textbf{In present method, these elements are seen as leaf nodes of binary tree and searching the binary tree in classical top-down way gives optimum run-time.} } 

\comment{
However computer implementation of the same was done with computational efficiency in mind and hence in actual values $\alpha^{2^id}$ for all $ 0\leq i\leq r$ was also calculated only once and similarly $\beta^{2^id}$ for all $ 0\leq i\leq r$ should  also calculated once, requiring space  of $ \mathcal{O}(r)$. However as indicated below last 3-4 levels of binary tree need not be traversed, but could be ascertained by considering all possible values at leaf nodes. 
\\
Basically square root of any $\alpha \in F^*_p $ will exist at leaf of binary tree of depth $(r-1)$ i.e. total possible values from which one is the correct value would be amongst the $2^{(r-1)}$ elements of the solution set.
\\
\\
$\sqrt{\alpha} \pmod p \in \{ \alpha^{(d+1)/2}\beta^{id} \text{   where   } 0\leq i\leq 2^{(r-1)}-1 \}$
\\
}

\subsection{How is it different from Tonelli-Shanks and AMM }
Let us start by analysing Tonelli-Shanks which starts with possible value of $\alpha^{d+1/2}$ as the possible root.\\
\textbf{Tonelli-Shanks }\\
\begin{algorithm}[H]
\caption{Tonell-Shanks Pseudo Code  to calculate $\sqrt{\alpha} \pmod p$. Returns 0 if none exits. Taken from wikipedia.}
\KwInput{$\alpha,\beta \in F_p^* , p$; $\beta $ is non-residue and $ p $ a prime.}
\KwOutput{$\sqrt{\alpha} $, if it exists else  $0$. }

\begin{algorithmic}
\IF {$(p=2)$}
	\RETURN {$\alpha$}
\ENDIF
\IF {$(\frac{\alpha}{p})= -1$} 
	\RETURN {$0$}
\ENDIF
\STATE $M \gets r , c \gets \beta^d , t \gets \alpha^d \text{ and }R \gets \alpha^{d+1/2} $ 
\WHILE {$t \not \equiv 1$}
	\STATE \colorbox{green} {$ \text{use repeated squaring to find minimum i}(0< i < M) \text{ such that } t^{2^i} \equiv 1 $}
	\STATE $b \gets c^{2^{M-i-1}}$
	\STATE $M \gets i , c \gets b^2 , t \gets tb^2 \text{ and } R \gets Rb$ 
\ENDWHILE
\RETURN{R}
\end{algorithmic}
\end{algorithm}
In Tonelli-Shanks algorithm the idea is to maintain, after each iteration,  $R,\alpha $  and  $t$  such that  $R^2\equiv \alpha t\pmod p$. $\alpha$ doesn't change after each iteration (since it's the quadratic residue whose root we're attempting to find), but  $R$  and  $t$  are updated. The heart of the algorithm lies in the fact that the order of  $t$  is strictly decreasing over iterations, so eventually it will have the order of  $1$. Only the element  $1 \pmod p$ is  order $1$. Once  $t \equiv 1\pmod p$, then we've found  $R^2 \equiv \alpha \pmod p$, where $R$ is our solution. Please note in best case scenario , order of $t$ decreases quickly and becomes $1$, however in worst case scenario order will decrease  by $1$ in each iteration. In each iteration one needs to find $i$ such that $t^{2^i}\equiv 1$ this is done by squaring $t$ till required  $i$ is found. This finding of $i$ at each iteration makes the Tonelli-Shanks algorithm's average run time, for certain non-residues , little worse than the proposed method. This difference can easily be felt for large values of $r$.\\
Please note that number of times outer `While loop' is executed depends on the  value of $\beta$. The execution time depend upon how fast order of $t$ decrease to 1. This is not only dependent on value of $\alpha$ but also on value of $\beta$.
\comment{ But even if order of t decreases 1 at a time  maximum number of times `While loop' executes is (r-1) .}
 \textbf{However it traverses the search tree bottom-up as is clear from  the nested loop coloured in green , inside the outer `While loop'}. Hence if chosen $\beta$ is such that number of times outer `While loop' is executed is around $r-1$ , the nested loop inside adds to  run time by $ \mathcal{O}{(r^2)}$.\\
\comment{
\\So basically, in every iteration of outer 'While loop', Tonelli-Shanks algorithm, is adjusting estimate ( initial estimate $= \alpha^{d+1/2}$) by scaling it with $\beta^{2^jd}$ where j is linear function of r and order of t. But once iteration is completed as it has revised value of t , it has to calculate the order of t so that it can calculate the value of new estimate $\beta^{2^jd}$ by which old estimate needs to be scaled.\\ 
\\This loop inside the while loop can be executed  maximum of r-1 times in one iteration . Hence in worst case this nested loop itself adds to $\mathcal{O}{(r^2)} $ run-time.} \textbf{This is what we mean, by saying,  Tonelli-shanks climbs the search tree and  traverses the search tree bottom up.  It needs to go up the tree ( with current value of $t$)  to the node where the $t^{2^i}\equiv 1$ , so that  it can estimate the value of $\beta^{2^jd } $ by which current estimate needs to be further scaled up.} If you don't have to climb the tree , means you have already reached the estimate. \\
For certain non residues $ \beta$, Tonelli-Shanks can do a better job, but proposed method proposes an algorithm which is  $ \beta$  agnostic and on average is likely to do a better job. In fact testing Tonelli-Shanks with pair ($\alpha,\beta$) when most of the time when square root of the node is -1,  illustrates the point very well. \\
In fact Tonelli-Shanks achieves optimum run-time or worst case run time depending upon pair  ($\alpha,\beta$) i.e. it depends upon the value of non-residue provided as input . \\

Our Proposed method basically traverses the binary tree top-down in deterministic fashion. There are at maximum $(r-1)$ steps involved to reduce $\alpha^{(p-1)/2} \equiv 1 $ to expression $\alpha^d* \beta^{2md} \equiv 1 $, with each step halving the exponent of $\alpha$. 
\comment {
Also proposed method starts by recognising through Theorem 2, the possible  solution set, one of whose  element must be square root of quadratic residue $\alpha$. So, in this case, algorithm has to basically search the set with maximum possible elements equal to $2^{(r-1)}$, efficiently. This is done easily and naturally  by a top down approach. 
In proposed method, however  for computational advantage, root of the tree is found only one time, starting from the leaf $\alpha^d$. Once search tree's root has been found, after that tree is traversed top down taking the branch on which side the desired leaf is going to be (depending upon if the square root is $1$ or $-1$). After first root of the search has been located (i.e max value of $i$ for which  $\alpha^{2^id} \equiv 1 $), it quickly traverses down to leaf node without going up and down the binary tree.\\

\\Computation advantage over Tonelli-Shanks occurs because tree is not traversed bottom up for each new estimate, rather it takes a definite path down, depending upon the square root of the node is $1$ or $-1$. However if $\alpha,\beta$ are such that square root of most nodes are $1$, proposed method would take same number of steps but  Tonelli-shanks will have an advantage ( these are best case scenario for Tonelli-Shanks) as it has incremented the guess correctly. But in scenario where nodes appear randomly as $1$ or $-1$ , number of steps in proposed method remain the same but run time of Tonelli-Shanks increases by factor of $r^2$. Hence, on average, Tonelli-Shanks performs worse than the proposed method. The contrast would be easily visible for large value of $r$. \textbf{In proposed  method we climb up  tree only for the first time and after that we just descend (no multiple climbs after adjusting the value )}. In fact, proposed method takes definitive way down from top to bottom leaf (to the solution), makes it possible to implement it much more efficiently. \textbf{ Tonelli-Shanks improves the result in each iteration by guessing in the right direction , whereas proposed method works towards the result in definitive way ; Hence proposed method  renders itself to more efficient implementation.} \\
}

\\Adleman, Manders \& Millers  \cite{AMM77}, is  the generalisation of Tonelli-Shanks for calculating $r^{th}$ root and it also tries to find out  $\alpha_j^{2^id} \equiv 1$ by each time revising the estimate  $\alpha_j$,  starting from the leaf till  it reaches the next root of the search tree. Its mechanism to  search the binary search tree to find out square root \cite{AMM77} is exactly similar to Tonelli-Shanks, hence also computationally expensive than the proposed method.\\
\\\textbf{The proposed method has the mechanism of moving up-down the search tree while, Tonelli-Shanks \cite {Tonelli} and  Adleman, Manders \& Millers \cite{AMM77} traverse the tree bottom up only.} Also proposed method is able to create a  pre calculated set (by Theorem 2)  which must contain root. Each decision at any intermediate node of the tree while travelling down the tree cuts the solution space by half.  Hence implementation of proposed method lends itself to much more efficient algorithm.

In fact insight given by the proposed algorithm allows one to calculate $\sqrt{\alpha}$ with out knowledge of strict non-residue. Availability of \emph{`relative'} non-residue to $\alpha$ would suffice for proposed method. Such  \emph{`relative'} non-residue can be found with probability of $\frac{2}{3} $ (averaged over all possible values of $\alpha $) compared to $\frac{1}{2} $.  Please see sections 2.4, 2.9 and 2.10 for  \emph{`relative'} non-residue and  example F in section 4.3.

\section{Run time analysis }
We can implement proposed algorithm in $ \mathcal{O}(r^2)$ and in extreme cases (when $r=\log{p}$ ) to  $ \mathcal{O}(\log^2{p})$ . This is achieved if  we are ready to store $r$ results of $(\alpha^d,\alpha^{2d}, \cdots, \alpha^{2^{(r-1)}d} ) $ and  $(\beta^d,\beta^{2d}, \cdots, \beta^{2^{(r-1)}d} )$ which we will need to calculate any way  to find out if $\alpha $ is indeed quadratic residue and $\beta$ is indeed a quadratic non residue. The storage requirement would be $ \mathcal{O}(r)$ and in worst case $ \mathcal{O}(\log{p})$ when  $p=2^r+1 $.  
The height of binary tree traversed would be at most $(r-1) $ and hence loop would be executed $(r-1)$ times. At each step in the loop, as values of $ \alpha^d, \alpha^{2d}, \cdots$ and  $ \beta^d, \beta^{2d}, \cdots$ are available, no more than $r$ multiplications would take place. Hence  run time of $ \mathcal{O}(r^2) $  which in worst case scenario would asymptotically approach $ \mathcal{O}(\log^2{p})$ i.e. when $r$ approaches $\log{p}$. Total number of multiplications required in worst case scenario are approx. $= 2*\log {p} + r^2$.\\
It is easy to see from theorem 2 proved in preceding section that  square root of any quadratic residue $\alpha$ is necessarily element   of the following set, if minimum $i$ for which $\alpha^{2^id} \equiv 1 $  and which say is  $i=r-j$. Please note least value for $j$ can be $1$.\\ 
$  \sqrt{\alpha} \in \{ (\alpha)^{(d+1)/2}\beta^{2^{(r-1-i)}kd} \}$ where $0 \leq k \leq 2^i-1 $; and $0\leq i \leq (r-1)$ \\ 
\\
This simplifies to \\
$  \sqrt{\alpha} \in \{(\alpha)^{(d+1)/2}\beta^{2^{(j-1)}kd} \}$  where $0 \leq k \leq 2^{(r-j)}-1 $; and $1 \leq j \leq r $ \\
\\
Hence  for practical purposes where $r-j$ is less than a suitable value (say 6) trying all possible  $2^{(r-j)} $ might be computationally cheaper. Hence for primes with smaller values of $r$  running time would be much faster with out any extra storage requirements at all.\\

\textbf{Also, the proposed algorithm is much more suitable to parallel computing as subsets of solution sets can be searched independently. This has not been practically tested though, linear improvement in run-time is expected with proportional increase in number of threads/CPUs.} 

\subsection{Comparison with Tonelli-Shanks}

Following is the table showing timings for checking first 10000 and 100000 elements for quadratic residue and it it exists calculating the same. 
Each result was verified. The proposed algorithm was implemented in python with recursive function to traverse the tree. Idea was to roughly compare the run time  with standard python code of tonelli-shanks algorithm available at Rosetta code (\url{http://rosettacode.org/wiki/Tonelli-Shanks_algorithm#Python}), primes (50 to 200 digits at \url{https://primes.utm.edu/lists/small/}). Both functions were given a non-residue as input too, so that comparison of basic algorithm can be made.\\Timings are in seconds = Tonelli-Shanks/ Our-proposed method \\
\\
1.  50  digit prime (till 10000 residues) = 1.54/ 1.24\\
2.  50  digit prime (till 100000 residues)= 15.01/12.36\\
3.  110 digit prime (till 10000 residues) = 7.18/5.69\\
4.  110 digit prime (till 100000 residues)= 72.93/62.23\\
5.  120 digit prime (till 10000 residues) = 8.85/6.91\\
6.  120 digit prime (till 10000 residues) = 88.79/77.39\\
7.  130 digit prime (till 100000 residues)= 9.59/7.6\\
8.  130 digit prime (till 10000 residues) = 88.71/75.88\\
9.  140 digit prime (till 10000 residues) = 13.05/10.43\\
10. 140 digit prime (till 100000 residues)= 128/106.1\\
11. 150 digit prime (till 10000 residues) = 12.76/10.41\\
12. 150 digit prime (till 100000 residues)= 137.76/117.81\\
13. 200 digit prime (till 5000 residues)  = 13.72/11.05\\
14. 200 digit prime (till 10000 residues) = 27.18/21.62\\
15. 200 digit prime (till 100000 residues)= 304.93/247.73\\
\\The improvement in time can be seen in all cases. This improvement was noticed without factoring the advantage of \emph{`relative'} non-residue, which proposed method enjoys. Improvement of $15-20\%$ in time was observed. However this advantage could be mainly because of efficient implementation also, although the  work done with in each iteration of loop in proposed algorithm seems to be less than the work done in each iteration of  Tonelli-Shanks. However more testing is required specially with larger value of $r$ (proth's prime) to see if the advantage vis a vis Tonelli-Shanks widen's and by how much.\\
\\It will not be out of place to mention here that for special type of large primes with very large $r$, Cipolla-Lehmer method performs better as it has asymptotically better run time. But the algorithm of Tonelli and Shanks for computing square roots modulo a prime number is the most used, and probably the fastest among the known algorithms when averaged over all prime numbers \cite {Tornaria}.
\subsection{ Examples}
Calculations have been shown, as would be easier to understand for humans. Exactly implementing this process may not give most efficient computer program, as it calculates things over and over again. However, process shows the simplicity of the basic algorithm.

\subsubsection{Example A : Let us try to calculate square root of 2 modulo 97 (given non-residue 5) }

First we check if 2 is indeed quadratic residue. we see $2^{48} \pmod{ 97} \equiv 1 $\\
We calculate  $2^{24} \pmod {97} \equiv  -1 $\\
As $2^{24} \pmod {97} \equiv  -1 $, multiplying both sides by identity $5^{48} \pmod {97} \equiv -1$\\
We get $2^{24*}5^{48} \pmod {97} \equiv 1 $\\
As exponent of 2 is still even again taking square root\\
We calculate $2^{12}*5^{24} \pmod {97} \equiv  -1 $\\
again multiplying both sides by identity $5^{48}\pmod {97} \equiv -1$\\
We get $2^{12}*5^{24}*5^{48} \pmod {97} \equiv 1 $\\
As exponent of 2 is still even again taking root\\
We get we get $2^{6}*5^{12}*5^{24} \pmod {97} \equiv -1 $\\
As $2^{6}*5^{12}*5^{24} \pmod {97} \equiv -1 $, multiplying both sides by identity $5^{48} \pmod {97}\equiv -1$\\
We get $2^{6}*5^{12}*5^{24}*5^{48} \pmod {97} \equiv 1 $\\
As exponent of 2 is still even again taking root\\
We calculate $2^{3}*5^{6}*5^{12}*5^{24} \pmod {97} \equiv -1 $\\
As the above expression is -1, multiplying both sides by identity $5^{48} \pmod {97}\equiv -1$\\
We get  $2^{3}*5^{6}*5^{12}*5^{24}*5^{48} \pmod {97} \equiv 1 $\\
As exponent of 2 is now odd we can multiply both sides by 2\\
We get $2^{4}*5^{6}*5^{12}*5^{24}*5^{48} \pmod {97} \equiv 2 $\\
Taking square root final time  we get \\
$\pm 2^{2}*5^{3}*5^{6}*5^{12}*5^{24} \pmod {97} \equiv \sqrt{2} $\\
\textbf{i.e$ \sqrt{2} \equiv \pm 83 $}\\
We find, indeed  $83*83 \pmod {97} \equiv 2$.\\

\subsubsection{Example B: Let us try to calculate square root of 6 modulo 43 (given non-residue 3)}

First we check if 6 is indeed quadratic residue. we see $6^{21} \pmod {43} \equiv 1 $\\
As the exponent of 6 is already odd, we need not take further square root. we multiply both sides by 6 \\
we get $ 6^{22} \pmod {43} \equiv 6 $\\
Taking square root for final time we get \\
implying $\sqrt{6} \equiv \pm 6^{11}$
\textbf{i.e. $\sqrt{6} \equiv \pm 36$}\\
We find, indeed  $36*36 \pmod {43} \equiv 6$\\
Please note we did not need to make use of quadratic non residue here at all.This is true for all primes of type $3 \pmod 4 $.\\
\subsubsection{Example C: Let us try to calculate square root of 2 modulo 41 (given non-residue 3)}

First we check if 2 is indeed quadratic residue. we see $2^{20} \pmod{41} \equiv 1 $\\
As the exponent of 2 is even we take square root of both sides\\
We get $2^{10} \pmod{41} \equiv -1 $\\
As result is -1 we multiply both sides by $3^{20} \equiv -1 $\\
We get $2^{10}* 3^{20} \pmod{41} \equiv 1 $\\
As the exponent of 2 is even we take square root of both sides\\
We get $2^{5}* 3^{10} \pmod{41} \equiv 1 $\\
As the exponent of 2 is odd now, we multiply both sides by 2\\
We get $2^{6}* 3^{10} \pmod{41}\equiv 2 $\\
Taking square root of both sides for final time \\
We get $\sqrt{2}  \equiv \pm 2^{3}* 3^{5} \pmod{41}  $\\
Implying  \textbf{$\sqrt{2}  \equiv \pm 17 \pmod{41}  $\\}
We find, indeed  $17*17 \pmod{41} \equiv 2$.\\

\subsubsection{Example D: Let us try to calculate square root of 12 modulo 13 (given non-residue 5)}

First we check if $12$ is indeed quadratic residue. we see $(12)^{6} \pmod{13} \equiv 1 $.\\
As the exponent of $12$ is even we take square roots of both sides 
We get $(12)^{3}\pmod{13} \equiv -1 $.\\
As the result is $-1$, we multiply both sides by $5^6\pmod{13} \equiv -1 $ \\
We get $(12)^{3}*5^{6} \pmod{13} \equiv 1 $.\\
As the exponent of $12$ is odd we multiply both sides by $12$\\
We get $(12)^{4}*5^{6} \pmod{13} \equiv 12 $.\\
Taking square root for the final time\\
We get $\sqrt{12} \equiv \pm (12)^{2}*5^{3} \pmod{13} $\\
implying $\sqrt{12} \equiv \pm  8 $\\
We find, indeed  $8*8 \pmod{13} \equiv 12 $\\
\\
\subsubsection{Example E: Let us try to calculate square root of 5 modulo $43^3=68921$ (given non-residue 3)}

We note $\phi(43^3)= 67240$ .We check and find  5 is indeed quadratic residue $5^{67240/2} = 5^{33620}  \equiv 1 \pmod{43^3}$ and 3 is indeed non residue $3^{33620}  \equiv -1 \ (mod\ 43^3). $ 
Now as exponent of 5 is even, we take square root \\
We get  $5^{16810}  \equiv -1 \pmod{43^3}$ 
As the result is -1, we multiply both sides by $3^{33620}  \equiv -1 \pmod{43^3} $ \\
we get $5^{16810}*3^{33620}  \equiv 1 \pmod{43^3}$\\
Now as exponent of 5 is even, we take square root \\
we get $5^{8405}*3^{16810}  \equiv -1 \pmod{43^3}$\\
As the result is -1, we multiply both sides by $3^{33620}  \equiv -1 \pmod{43^3} $ \\
We get $5^{8405}*3^{16810}*3^{33620}  \equiv 1 \pmod{43^3}$\\
Now as exponent of 5 is odd we multiply both sides with 5 and take square root for the last time \\
We get $ \sqrt{5} \pmod{43^3} \equiv 5^{8405}*3^{8405}*3^{16810} $\\
Implying \textbf{$ \sqrt{5} \pmod{43^3} \equiv \pm 3226$}\\
We find,  indeed  $3226*3226 \pmod{43^3} \equiv 5 $.\\

\subsubsection{Example F: Let us try to calculate square root of 6 modulo 97 (No non residue given)}

First we check if 6 is indeed quadratic residue. we see $6^{48} \pmod {97} \equiv 1 $.\\
As the exponent of 6 is even we take square roots of both sides \\
We get $6^{24} \pmod {97} \equiv 1 $\\
As the exponent of  6 is even we take square roots of both sides \\
We get $6^{12} \pmod {97} \equiv 1 $\\
As the exponent of 6 is even we take square roots of both sides \\
We get $6^{6} \pmod {97} \equiv -1 $\\
Now we do not have any quadratic non residue given, we select any random number hoping it would be quadratic non residue (We note probability of finding a non-residue is $\frac{1}{2} $).  We will deal it in different parts. \\

\paragraph{Part 1 : (We get 2 as randomly selected number)\\ }
Although 2 is non residue but we see $f(2) > f(6) $ (as $f(2)=3$ and $f(6)=1$), hence 2 is \emph{`relative'} non-residue and  suitable. Also $f(2)=3 \text{ implies } 2^{24} \equiv -1 $.\\
As the value $6^{6} \pmod {97} \equiv -1 $, we multiply both sides by $ 2^{24} \equiv -1 $ \\
We get $6^{6} * 2^{24} \pmod {97} \equiv 1 $\\
As the exponent 6 is even we take square roots of both sides \\
We get $6^{3}*2^{12}\pmod {97} \equiv -1 $\\
As the value is $-1$, we multiply both sides by $2^{24} \equiv -1$ \\
we get $6^{3}* 2^{12} *2^{24} \pmod {97} \equiv 1 $\\
As exponent of 6 is now odd we can multiply both sides by 6\\
we get $6^{4}* 2^{12} *2^{24} \pmod {97} \equiv 6 $\\
Taking square root for the final time, 
we get $6^{2} *2^{6} *2^{12}\pmod {97} \equiv \sqrt(6) $\\
Implying  \textbf{$\sqrt{6}  \equiv \pm 54 \pmod {97}  $\\}
We find,  indeed  $54*54 \pmod {97} \equiv 6$.
\paragraph{Part 2 : (We get 9 as randomly selected number)\\ }
Although 9 is non residue but we see $f(9) > f(6) $ (as $f(9)=2$ and $f(6)=1$), hence 9 is \emph{`relative'} non-residue and  suitable. Also $f(9)=2 \text{ implies } 9^{12} \equiv -1 $.\\
As the value $6^{6} \pmod {97} \equiv -1 $, we multiply both sides by $ 9^{12} \equiv -1 $ \\
We get $6^{6} * 9^{12} \pmod {97} \equiv 1 $\\
As the exponent 6 is even we take square roots of both sides \\
We get $6^{3} * 9^{6} \pmod {97} \equiv 1 $\\
As exponent of 6 is  now odd we can multiply both sides by 6\\
We get $6^{4} * 9^{6} \pmod {97} \equiv 6 $\\
Taking square root for the final time, 
We get $6^{2} * 9^{3}\pmod {97} \equiv \sqrt(6) $\\
Implying  \textbf{$\sqrt{6}  \equiv \pm 54 \pmod {97}  $\\}
We find,  indeed  $54*54 \pmod {97}\equiv 6$.
\paragraph{Part 3 : (We get 22 as randomly selected number)\\ }
Although 22 is also non-residue but we see $f(22) = f(6) $ (as $f(22)=1$ and $f(6)=1$), hence 22 is not suitable for calculating square root of 6.  We will have to choose another number and see if it fits the bill!\\ 
\\The basic idea was to show that probability of finding a suitable \emph{`relative'} non-residue  for calculating  the  square root of a given number, is greater ( by around 17\% , averaged over all residues )   than finding non-residues, as is required in other standard methods. \\

\section{ Future ideas }
The strength of this algorithm is its simplicity and general applicability to all odd primes. However main drawback of this algorithm is the presupposition of the quadratic non-residue. But this can be over come to an extent. This algorithm does not necessarily demands a quadratic non residue. In fact given any element $ \alpha$  whose square root is to be calculated with minimum $i_\alpha $  for which  $\alpha^{2^{i_\alpha}d} \equiv -1 $, availability of any  $\gamma$ (not necessarily a quadratic non residue) with minimum $i_\gamma $ for which  $\gamma^{2^{i_\gamma}d} \equiv -1 $ such that $ i_\gamma > i_\alpha$ will suffice. It is easy to see that if $\gamma$ is indeed non residue, this condition will always hold for any quadratic residue $\alpha$. Please see example 'F' in preceding section\\
\\In fact depending upon the height of tree being examined availability of an element whose forth root or for that matter  eighth root does not exist (although square root may exist, i.e. it is a quadratic residue ) should also suffice(depending upon height of tree to be explored. We can easily see such \emph{`relative'} non-residues should be far easier to locate, as of all the elements 3/4th of all elements would be fourth root non residue and 7/8th of all elements would be eighth root non residue. It has been shown in one of the preceding section,  \textbf{`Discussion about \emph{`relative'} non-residue} probability of finding such \emph{`relative'} non-residue, averaged over all residues, increases from $\frac{1}{2} $ to $ \frac{2}{3}$.\\
\\However more fundamental area would be to find a way to avoid using quadratic non residue altogether maintaining  similar kind of run-time performance. \textbf{From philosophical point of view, mandatory requirement of an element whose square root can not be calculated (i.e. non-residue), for finding square root of all other numbers (quadratic residues), definitely begs a more fundamental answer.} \\
\\Another area of exploration could be to quickly locate the desired solution given the solution set. For a given quadratic residue $\alpha$  whose square root is to be calculated (and also given a non-residue $\beta$), we know the set of possible solutions(cardinality at most equal to $2^{r-1}$). Challenge is to find an efficient search mechanism to quickly locate the desired result. Although it seems like discrete logarithm problem at first glance, but a serious look may give better search time algorithms.\\

\section{Conclusion}

Ability to find  square roots efficiently, in a finite field has its applications in cryptosystem  as broader classes  of elliptic and hyper-elliptic curve cryptosystem can be set up more efficiently.
Another idea could be to use it for efficient deterministic primality testing for Proth primes as shown by Sze \cite{Sze}.\\
\comment{
\\I would like to thank Prof. Ritabarata Munshi, Prof. Palash Sarkar,  Prof Neena Gupta of ISI Kolkata, Prof. G.P. Biswas of IIT Dhanbad,  Shreesh Maharaj of Swami Vivekanand University, Belur math and  Malay Khandelwal of thinkC  for not only  initiating me to number theory  but also for  motivation, guidance and encouragement. Special thanks are due to my colleagues Gaurav Sinha, IRS and Syed Waquar Raza, IPS for discussing the paper, taking pains to read and re-read the paper endless number of times and suggesting valuable improvements.\\}

\end{document}